\magnification=\magstep1
\hsize=16.5 true cm 
\vsize=23.6 true cm
\font\bff=cmbx10 scaled \magstep1
\font\bfff=cmbx10 scaled \magstep2
\font\bffg=cmbx10 scaled \magstep3

\font\smc=cmcsc10 
\parindent0cm
\def\cl{\centerline}           %
\def\rl{\rightline}            %
\def\bp{\bigskip}              %
\def\mp{\medskip}              %
\def\sp{\smallskip}            %
\def\bc{{\bf c}}
           %
\def\Bbb#1{\hbox{\boldmas #1}} %
\def\Q{\Bbb Q}                 %
\def\R{\Bbb R}                 %
\def\N{\Bbb N}                 %
\def\Z{\Bbb Z}                 %
\expandafter\edef\csname amssym.def
\endcsname{%
       \catcode`\noexpand\@=\the\catcode`\@\space}

\catcode`\@=11

\def\undefine#1{\let#1\undefined}
\def\newsymbol#1#2#3#4#5{\let\next@\relax
 \ifnum#2=\@ne\let\next@\msafam@\else
 \ifnum#2=\tw@\let\next@\msbfam@\fi\fi
 \mathchardef#1="#3\next@#4#5}
\def\mathhexbox@#1#2#3{\relax
 \ifmmode\mathpalette{}{\m@th\mathchar"#1#2#3}%
 \else\leavevmode\hbox{$\m@th\mathchar"#1#2#3$}\fi}
\def\hexnumber@#1{\ifcase#1 0\or 1\or 2\or 3\or 4\or 5\or 6\or 7\or 8\or
 9\or A\or B\or C\or D\or E\or F\fi}

\font\tenmsa=msam10
\font\sevenmsa=msam7
\font\fivemsa=msam5
\newfam\msafam
\textfont\msafam=\tenmsa
\scriptfont\msafam=\sevenmsa
\scriptscriptfont\msafam=\fivemsa
\edef\msafam@{\hexnumber@\msafam}
\mathchardef\dabar@"0\msafam@39
\def\dashrightarrow{\mathrel{\dabar@\dabar@\mathchar"0\msafam@4B}}
\def\dashleftarrow{\mathrel{\mathchar"0\msafam@4C\dabar@\dabar@}}

\def\ulcorner{\delimiter"4\msafam@70\msafam@70 }
\def\urcorner{\delimiter"5\msafam@71\msafam@71 }
\def\llcorner{\delimiter"4\msafam@78\msafam@78 }
\def\lrcorner{\delimiter"5\msafam@79\msafam@79 }
\def\yen{{\mathhexbox@\msafam@55}}
\def\checkmark{{\mathhexbox@\msafam@58}}
\def\circledR{{\mathhexbox@\msafam@72}}
\def\maltese{{\mathhexbox@\msafam@7A}}

\font\tenmsb=msbm10
\font\sevenmsb=msbm7
\font\fivemsb=msbm5
\newfam\msbfam
\textfont\msbfam=\tenmsb
\scriptfont\msbfam=\sevenmsb
\scriptscriptfont\msbfam=\fivemsb
\edef\msbfam@{\hexnumber@\msbfam}
\def\Bbb#1{{\fam\msbfam\relax#1}}
\def\widehat#1{\setbox\z@\hbox{$\m@th#1$}%
 \ifdim\wd\z@>\tw@ em\mathaccent"0\msbfam@5B{#1}%
 \else\mathaccent"0362{#1}\fi}

\def\widetilde#1{\setbox\z@\hbox{$\m@th#1$}%
 \ifdim\wd\z@>\tw@ em\mathaccent"0\msbfam@5D{#1}%
 \else\mathaccent"0365{#1}\fi}
\font\teneufm=eufm10
\font\seveneufm=eufm7
\font\fiveeufm=eufm5
\newfam\eufmfam
\textfont\eufmfam=\teneufm
\scriptfont\eufmfam=\seveneufm
\scriptscriptfont\eufmfam=\fiveeufm

\newsymbol\risingdotseq 133A
\newsymbol\fallingdotseq 133B
\newsymbol\complement 107B
\newsymbol\nmid 232D
\newsymbol\rtimes 226F
\newsymbol\thicksim 2373

\font\eightmsb=msbm8   \font\sixmsb=msbm6   \font\fivemsb=msbm5
\font\eighteufm=eufm8  \font\sixeufm=eufm6  \font\fiveeufm=eufm5
\font\eightrm=cmr8     \font\sixrm=cmr6     \font\fiverm=cmr5
\font\eightbf=cmbx8    \font\sixbf=cmbx6    
      \font\eighti=cmmi8   \font\sixi=cmmi6
\font\ninesy=cmsy9     \font\eightsy=cmsy8  \font\sixsy=cmsy6
     \font\eightit=cmti8  
     \font\eightsl=cmsl8  
     \font\eighttt=cmtt8

\font\eightsmc=cmcsc8
\newskip\ttglue
\newfam\smcfam
\def\eightpoint{\def\rm{\fam0\eightrm}%
  \textfont0=\eightrm \scriptfont0=\sixrm \scriptscriptfont0=\fiverm
  \textfont1=\eighti \scriptfont1=\sixi \scriptscriptfont1=\fivei
  \textfont2=\eightsy \scriptfont2=\sixsy \scriptscriptfont2=\fivesy
  \textfont3=\tenex \scriptfont3=\tenex \scriptscriptfont3=\tenex
  \def\smc{\fam\smcfam\eightsmc}
  \textfont\smcfam=\eightsmc          
\textfont\eufmfam=\eighteufm              \scriptfont\eufmfam=\sixeufm
     \scriptscriptfont\eufmfam=\fiveeufm
\textfont\msbfam=\eightmsb            \scriptfont\msbfam=\sixmsb
     \scriptscriptfont\msbfam=\fivemsb
\def\it{\fam\itfam\eightit}%
  \textfont\itfam=\eightit
  \def\sl{\fam\slfam\eightsl}%
  \textfont\slfam=\eightsl
  \def\bf{\fam\bffam\eightbf}%
  \textfont\bffam=\eightbf \scriptfont\bffam=\sixbf
   \scriptscriptfont\bffam=\fivebf
  \def\tt{\fam\ttfam\eighttt}%
  \textfont\ttfam=\eighttt
  \tt \ttglue=.5em plus.25em minus.15em
  \normalbaselineskip=9pt
  \def\MF{{\manual opqr}\-{\manual stuq}}%
  \let\big=\eightbig
  \setbox\strutbox=\hbox{\vrule height7pt depth2pt width\z@}%
  \normalbaselines\rm}
\def\eightbig#1{{\hbox{$\textfont0=\ninerm\textfont2=\ninesy
  \left#1\vbox to6.5pt{}\right.\n@space$}}}

\catcode`@=13 
\cl{\bffg On real functions with graphs either}
\sp
\cl{\bffg connected or locally connected}
\bp
\cl{\bfff Gerald Kuba}
\bp\mp
\vbox{\eightpoint {\bf Abstract.} Let $\,{\cal G}\,$ denote the 
family of all subspaces $\,G\,$ of the plane $\,\R^2\,$
such that $\,G\,$ is the graph of a function from $\R$ to $\R$.
We prove that $\,{\cal G}\,$ 
has two subfamilies $\,{\cal G}_1,{\cal G}_2\,$ of {\it connected}
spaces such that the cardinality of $\,{\cal G}_1\,$ is $\,\bc:=2^{\aleph_0}\,$ 
and the cardinality of $\,{\cal G}_2\,$ is $\,2^\bc\,$, 
every space in $\,{\cal G}_1\,$ is {\it completely metrizable},
each $\,G\in{\cal G}_2\,$ is a dense subset of $\,\R^2$,
and if $\,X_1,X_2\,\in\,{\cal G}_1\cup{\cal G}_2\,$ are distinct
then the space $\,X_1\,$ is neither homeomorphic to a subspace
of $\,X_2\,$ nor homeomorphic to a proper subspace of $\,X_1\,$.
On the other hand, the family $\,{\cal G}\,$ contains precisely 
$\,\aleph_0\,$ {\it locally connected} spaces
up to homeomorphism, and if $\,X,Y\,$ are such spaces (including 
the case $X=Y$) then 
$\,X\,$ is homeomorphic to some proper subspace of $\,Y\,$.
Furthermore, if $\,\tau\,$ is a topology on the set $\,\R\,$
finer than the Euclidean topology and the space
$\,(\R,\tau)\,$ is separable and locally connected
then the space is locally compact and 
homeomorphic to some space in $\,{\cal G}\,$.
In a very natural way we 
establish a complete classification of all these refinements $\,\tau\,$
of the real line. 
\sp 
}
\bp\mp
{\bff 1.~Introduction}
\mp
Let $\,|S|\,$ denote the cardinal number (the {\it size}) of any set $\,S\,$
and let $\,\bc:=|\R|=2^{\aleph_0}\,$ denote the cardinality of the continuum.
A {\it real function} is a function from $\,\R\,$ to $\,\R\,$.
Let $\,{\cal R}\,$ denote 
the family of all real functions. Naturally, $\,|{\cal R}|=2^\bc\,$.
As usual, any function is 
regarded as identical with its graph, whence $\,F\,$ is a {\it subset} 
of $\,\R^2\,$ for every $\,F\in{\cal R}\,$.
In the following, we
consider $\,F\in{\cal R}\,$ as a {\it subspace} of the Euclidean plane
$\,\R^2\,$. 
In doing so every $\,F\in{\cal R}\,$
induces an interesting topology on the real line, see Section 2.
\sp
It is evident that a real function is continuous 
if and only if it is a pathwise connected subset of the plane.
While it is clear that all pathwise connected spaces $\,F\in{\cal R}\,$
are homeomorphic, this is not true for 
connected spaces $\,F\in{\cal R}\,$. 
Actually, our first goal is to prove the following theorem
about connected real functions.
For abbreviation, let us call topological 
spaces $\;X_i\;(i\in I)\;$ {\it incomparable} if and only if 
for $\,i,j\in I\,$ the space 
$\,X_i\,$ is homeomorphic to a subspace $\,S_j\,$ of $\,X_j\,$
only in the trivial case $\,X_i=S_j=X_j\,$.
In other words, spaces are incomparable when they are 
mutually non-homeomorphic 
and no space is homeomorphic to a proper subspace of itself or of any 
other space. 
\mp
{\bf Theorem 1.} {\it There exist $\,2^\bc\,$ real functions which are 
incomparable, connected, dense subspaces of $\,\R^2\,$.} 
\mp
It is plain 
that if a real function  $\,F\in{\cal R}\,$ is a {\it dense} 
subset of $\,\R^2\,$ (as in Theorem 1) then the function $\,F\,$
is {\it everywhere discontinuous.}
It is well-known (see [1]) that a subspace $\,X\,$ of $\,\R^2\,$
is {\it completely metrizable}
if and only if $\,X\,$ is a G$_\delta$-subset of $\,\R^2\,$.
Since it is clear that $\,\R^2\,$ has precisely $\,\bc\,$ 
G$_\delta$-subsets, at most $\,\bc\,$ 
functions provided by Theorem 1 are completely metrizable spaces. 
In fact, no space provided by Theorem 1 is completely metrizable
due to the following theorem.
\mp
{\bf Theorem 2.} {\it If $\,F\in {\cal R}\,$ is a connected, 
completely metrizable space and $\,S\,$ is the set of all 
$\,x\in\R\,$ such that the function $\,F\,$ is not continuous 
at the point $\,x\,$ then $\,S\,$ is a meager subset of $\,\R\,$.}
\mp
Note that the set $\,S\,$ in Theorem 2 is a F$_\sigma$-subset of $\,\R\,$,
see [7] 7.1. In view of Theorem 1 and Theorem 2 there 
arises a natural question which is 
answered by the following theorem.
\mp
{\bf Theorem 3.} {\it There exist
$\,\bc\,$ real functions which are 
incomparable, connected, completely metrizable spaces.} 
\bp
{\bff 2.~Refinements of the real line}
\mp
If $\,\tau\,$ is the topology of a space and
$\,x\,$ is a point then let $\,{\cal U}_\tau(x)\,$                   
denote the filter of all neighborhoods of $\,x\,$.
Let $\,\eta\,$ denote the natural Euclidean topology on 
the real line $\,\R\,$. Let $\,{\cal T}\,$ be the 
family of all topologies $\,\tau\,$ on $\,\R\,$ which are finer 
than $\,\eta\,$ or, equivalently, where 
$\,{\cal U}_\eta(x)\subset{\cal U}_\tau(x)\,\,$ for every $\,x\in\R\,$.
(Naturally, $\,(\R,\tau)\,$ is a Hausdorff space for every 
$\,\tau\in{\cal T}\,$.) 
For $\,\tau\in {\cal T}\,$ put 
$\,\Gamma(\tau)\,:=\,
\{\,x\in\R\;|\;{\cal U}_\eta(x)\not={\cal U}_\tau(x)\,\}\,$.
So we can say that {\it $\,\tau\,$ is strictly finer than $\,\eta\,$
at $\,x\,$} if and only if $\,x\in\Gamma(\tau)\,$. 
\mp
By considering the functions $\,F\in{\cal R}\,$ as 
subspaces of the Euclidean plane, in a very natural way
we obtain topologies in the family $\,{\cal T}\,$.
For $\,F\in{\cal R}\,$ define a topology $\,\tau[F]\,$ on $\,\R\,$
by declaring $\,U\subset\R\,$ open
if and only if $\;U\,=\,\{\,x\in\R\;|\;(x,F(x))\in V\,\}\;$
for some open subset $\,V\,$ of the plane $\,\R^2\,$.
In other words, $\,\tau[F]\,$ is the unique topology $\,\tau\,$ on $\,\R\,$
such that $\;x\mapsto (x,F(x))\;$ is a homeomorphism
from the space $\,(\R,\tau)\,$ onto the subspace $\,F\,$ of $\,\R^2\,$.
In particular, the space $\,(\R,\tau[F])\,$ is separable and
metrizable for every $\,F\in{\cal R}\,$. Furthermore, 
$\,(\R,\tau[F])\,$ is connected if and only if $\,F\,$ is connected,
and $\,(\R,\tau[F])\,$ is completely metrizable 
if and only if $\,F\,$ is a G$_\delta$-subset of $\,\R^2\,$.
\mp
A proof of the following important observation 
is straightforward.
\sp
(2.1)$\;${\it If $\,F\in{\cal R}\,$ and $\,\xi\in\R\,$ 
then $\;\xi\in\Gamma(\tau[F])\;$ if and only if 
$\,F\,$ is not continuous at $\,\xi\,$.}
\mp
In particular, $\,\tau[F]=\eta\,$ if and only if the real function $\,F\,$
is continuous. In view of (2.1), from Theorem 1 and Theorem 3 we derive 
the following two corollaries.
\mp
{\bf Corollary 1.} {\it There exist $\,2^\bc\,$ topologies 
$\,\tau\in{\cal T}\,$ with $\,\Gamma(\tau)=\R\,$ 
such that the corresponding spaces $\,(\R,\tau)\,$ 
are incomparable, connected, separable, and metrizable.}
\mp
{\bf Corollary 2.} {\it There exist $\,\bc\,$ topologies 
$\,\tau\in{\cal T}\,$ such that the corresponding spaces $\,(\R,\tau)\,$ 
are incomparable, connected, separable, and completely metrizable.}
\mp
The cardinality $\,2^\bc\,$ in Corollary 1 is maximal 
since a metric on $\,\R\,$ is a function from 
$\,\R\times\R\,$ into $\,[0,\infty[\;$ and, trivially,
there are precisely $\,2^\bc\,$ mappings from $\,X\,$ to $\,Y\,$ 
whenever $\,|X|=|Y|=\bc\,$.
Also the cardinality $\,\bc\,$ in Corollary 2 is maximal 
because incomparable spaces are mutually non-homeomorphic and
there exist precisely $\,\bc\,$ separable, completely metrizable spaces
up to homeomorphism (see [4]). 
\mp
{\it Remark.} By (2.1) and [7] 7.1, for every $\,F\in{\cal R}\,$
the set $\,\Gamma(\tau[F])\,$ is a F$_\sigma$-set of reals.
More generally, it is a nice exercise to verify that  
$\,\Gamma(\tau)\,$ is a F$_\sigma$-set for every {\it metrizable}
topology $\,\tau\in{\cal T}\,$. This is not true for an 
arbitrary topology $\,\tau\in{\cal T}\,$.
Actually, for {\it every} set $\,S\subset\R\,$
we can easily define a topology $\,\tau_S\in{\cal T}\,$ with 
$\;\Gamma(\tau_S)=S\,$ by declaring a set 
$\,U\,$ $\tau_S$-open if and only if $\;U=V\cup T\;$
for some Euclidean open $\,V\subset\R\,$ and some $\,T\subset S\,$. 
(The space $\,(\R,\tau_S)\,$ is obviously paracompact.
In view of the metrization theorem [1] 4.4.7 
it is evident that the space $\,(\R,\tau_S)\,$ is 
metrizable if $\,S\,$ is a F$_\sigma$-set.)
\mp\sp
Clearly, if $\,\tau\in{\cal T}\,$ and the space $\,(\R,\tau)\,$ is {\it not 
separable} or {\it not metrizable}
then $\,\tau=\tau[F]\,$ is impossible for any $\,F\in{\cal R}\,$.
A natural example of a topology $\,\tau\in{\cal T}\,$ which is 
metrizable but not separable is the discrete topology. 
In Section 9 we will apply Theorem 1 in order to construct 
extremely many incomparable, separable but not metrizable,
connected refinements of the real line. 
On the other hand, in the following we investigate 
an important subfamily $\,{\cal L}\,$
of $\,{\cal T}\,$ where every topology  $\,\tau\in{\cal L}\,$
coincides with $\,\tau[F]\,$ for some $\,F\in{\cal R}\,$. 
This family $\,{\cal L}\,$ is the family of all {\it locally connected}
topologies in $\,{\cal T}\,$.
\mp
Naturally, the property {\it locally connected} is quite different 
from the property {\it connected}. This difference 
is extreme in the realm of real functions or refinements of the real line,
respectively. Indeed, the following proposition shows that 
a discontinuous real function and, more generally, 
a strict refinement of the real line can never be both 
connected and locally connected. 

\mp
{\bf Proposition 1.} {\it If $\,\tau\in{\cal T}\,$ and
and the space 
$\,(\R,\tau)\,$ is connected and locally connected then $\,\tau\,$ 
is the Euclidean topology $\,\eta\,$.
Moreover, if $\,\xi\in\R\,$ and 
$\,\tau\in{\cal T}\,$ is strictly finer than $\,\eta\,$ at $\,\xi\,$ 
and the space $\,(\R,\tau)\,$ is connected
then the point $\,\xi\,$ does not have a local basis 
of connected sets in the space $\,(\R,\tau)\,$.}
\mp
{\it Proof.} Let $\,a\in\R\,$ and let 
$\,\tau\in{\cal T}\,$ and assume 
that $\,{\cal U}_\eta(a)\not={\cal U}_\tau(a)\,$. 
Let $\,{\cal B}_a\,$
be a local basis at $\,a\,$ which consists of 
$\tau$-connected, $\tau$-open sets. A fortiori, any $\tau$-connected 
set is $\eta$-connected and hence an interval.
Thus $\,{\cal B}_a\,$ contains an interval which is not $\eta$-open.
If, for example, $\,[x,y[\;$
is such an interval, then $\;[x,\infty[\,=[x,y[\,\cup\,]x,\infty[\;$
and $\;]{-\infty,x}[\;$ are disjoint $\tau$-open sets whose union equals 
$\,\R\,$. If $\,[x,y]\,$ with $\,x<y\,$ or $\,x=y\,$ is such an interval then 
$\,[x,y]\,$ is $\eta$-closed and hence both $\tau$-open and $\tau$-closed.
Generally, in any case we can be sure that 
the space $\,(\R,\tau)\,$ is not connected, {\it q.e.d.}
\mp
{\it Remark.} In view of (2.1) and by virtue of Proposition 1, 
all spaces depicted in Theorem 1 and in 
Corollary 1 are {\it nowhere locally connected.}
\mp
The difference between the properties 
{\it connected} and {\it locally connected} in the 
realm $\,{\cal R}\,$ 
or in the realm $\,\{\,(\R,\tau)\;|\;\tau\in{\cal T}\,\}\,$
can also be quantified by transfinite cardinal numbers 
when we compare Theorems 1 and 3 and Corollaries 1 and 2
with the following theorem.
\mp
{\bf Theorem 4.} {\it There are only countably many mutually non-homeomorphic 
topologies $\,\tau\,$ in $\,{\cal T}\,$ such that the space 
$\,(\R,\tau)\,$ is locally connected. If such a space 
$\,(\R,\tau)\,$ is separable then 
$\,\tau=\tau[F]\,$ for some real function $\,F\,$
and the set $\,\Gamma(\tau)\,$ is 
countable and Euclidean closed.}
\mp
It is very easy to find for each countably infinite 
closed set $\,S\subset\R\,$ a sequence $\;F_1,F_2,F_3,...\;$
of mutually non-homeomorphic 
locally connected spaces in $\,{\cal R}\,$ with $\,\Gamma(\tau[F_k])=S\,$
for every $\,k\in\N\,$.
Indeed, since $\,S\,$ is scattered, we can define 
mutually disjoint compact intervals $\;[a_k,b_k]\;(k\in\N)\;$
with $\,a_k<b_k\,$ and $\;S\cap [a_k,b_k]=\{a_k,b_k\}\,$.
Choose a bijection $\,\varphi\,$ from $\,S\,$ onto $\;\N\setminus\{1\}\;$  
and define for every $\,k\in\N\,$ a function 
$\,F_k\in{\cal R}\,$ by $\;F_k(x)=1\;$ if
$\;x\in[a_1,b_1]\cup\cdots\cup[a_k,b_k]\;$ 
and $\;F_k(x)=\varphi(x)\;$ if 
$\;x\,\in\,S\setminus\{\,a_1,b_1,...,a_k,b_k\}\;$
and $\;F_k(x)=0\;$ if 
$\;x\,\in\,\R\setminus(S\cup[a_1,b_1]\cup\cdots\cup[a_k,b_k])\,$.
While the spaces $\,F_n\,$ and $\,F_m\,$ are never 
homeomorphic for distinct $\,n,m\in\N\,$ (because  
the space $\,F_k\,$ has precisely $\,k\,$ 
infinite compact components), 
the spaces $\;F_k\;(k\in\N)\;$ 
are not incomparable. On the contrary, it is plain that
$\,F_n\,$ is homeomorphic to a proper subspace of $\,F_m\,$ 
whenever $\,n,m\in\N\,$. 
(If $\,n\geq m\,$ then embed the subspace 
$\,([a_1,b_1]\cup\cdots\cup[a_n,b_n])\times\{1\}\,$ 
of $\,F_n\,$ into one component $\,[a,b]\times\{1\}\,$
of $\,F_m\,$.)
Actually, locally connected spaces $\,F\in{\cal R}\,$ are never incomparable.  
Moreover, in Section~8 we will prove that 
every locally connected refinement of $\,\R\,$ is homeomorphic to 
some proper subspace of itself and that if 
$\,{\cal F}\,$ is a family of 
locally connected spaces $\,(\R,\tau)\,$ with $\,\tau\in{\cal T}\,$
such that $\,|{\cal F}|\geq 2\,$ and 
no $\,X\in{\cal F}\,$ is homeomorphic to a subspace
of any $\;Y\,\in\,{\cal F}\setminus\{X\}\;$
then $\,|{\cal F}|=2\,$ and one space in $\,{\cal F}\,$ is discrete 
and the other space in $\,{\cal F}\,$ is separable.  
\bp
{\bff 3.~Massive connected functions}
\mp
Naturally, a real function $\,F\,$ is a Lebesgue measurable
subset of $\,\R^2\,$ if and only if $\,F\,$ is a null set in $\,\R^2\,$.
In particular, if $\,F\,$ is a {\it Borel subset} of $\,\R^2\,$
then $\,F\,$ is null. 
Of course, all spaces $\,F\,$
provided by Theorem 3 are Borel sets. 
Since there are only $\,\bc\,$ Borel subsets of $\,\R^2\,$,
it is impossible that all spaces depicted in Theorem 1 
are Borel sets. We will prove Theorem 1 by constructing 
functions which all are far from being Lebesgue measurable.
\sp
Following [3], a set $\,A\subset\R^2\,$ is {\it massive}
if and only if $\;A\cap L\not=\emptyset\;$ for every set $\,L\subset\R^2\,$
of positive Lebesgue measure
(or, equivalently, for every non-null Lebesgue measurable set 
$\,L\subset\R^2\,$). In particular, every massive subset 
of $\,\R^2\,$ is {\it dense} in $\,\R^2\,$.
Clearly, if $\,A\subset\R^2\,$ is Lebesgue measurable and massive 
then $\;\R^2\setminus A\;$ is a null set.
Consequently, if $\,F\in{\cal R}\,$ is massive
then $\,F\,$ is not Lebesgue measurable.
Our goal is to prove a stronger version of 
Theorem 1 where the property {\it dense} is replaced 
with the property {\it massive}. 
\mp
A natural way to accomplish this is to combine the main ideas of
three different proofs 
using transfinite induction, namely the proofs of 
[6, \S 35, V.~Theorem 2] and 
[3, Ch.~8, Theorem 5] and [5, Theorem 1].
Since this has to be done in a bit tricky way,
for the sake of better reading we conclude this section with 
four lemmas pivotal for the considerations which  
prove our sharper version of Theorem~1. In the next section we carry out 
the remaining part of the proof of Theorem 1 in detail.
In the following let $\,\pi\,$ denote the projection $\;(x,y)\mapsto x\;$
from $\,\R^2\,$ onto $\,\R\,$.
\mp
{\bf Lemma 1.} {\it If $\,B\subset\R^2\,$ is a Borel set which is not null
in $\,\R^2\,$ then $\,|\pi(B)|=\bc\,$.}
\mp
{\it Proof.} 
Let $\,B\subset\R^2\,$ be a Borel set.
Assume that $\,|\pi(B)|<\bc\,$. Then $\,\pi(B)\,$ must be countable 
because $\,\pi(B)\,$ is an {\it analytic} set
and (see [2] 11.20) every uncountable analytic set is of size $\,\bc\,$.
Consequently, $\,\pi(B)\times\R\,$ is a countable union of straight lines 
and hence a null set in the plane. 
Since $\,B\,$ is a subset of $\,\pi(B)\times\R\,$, we conclude that
$\,B\,$ is null, {\it q.e.d.} 
\mp\sp
{\bf Lemma 2.} {\it If $\,F\in{\cal R}\,$ and $\,A\cap F\not=\emptyset\,$
whenever $\,A\,$ is a closed subset of $\,\R^2\,$ and $\,\pi(A)\,$
contains a nondegenerate interval
then $\,F\,$ is connected.}
\mp
{\it Proof.} Assume indirectly that $\,F\,$ is not connected.
Then we can find open sets $\;U,V\subset\R^{2}\;$ such that 
$\,U\cap F\,$ and $\,V\cap F\,$ are nonempty and disjoint and 
$\;U\cup V\supset F\,$. Consequently, $\;U\cap V=\emptyset\;$
since $\,F\,$ is a dense subset of $\,\R^{2}\,$.
Since the projection $\,\pi\,$ is an open mapping,
the two sets $\;U':=\pi(U)\;$ and $\,V':=\pi(V)\;$ 
are open subsets of $\,\R\,$. Since $\;F \subset U\cup V\;$
and $\,F\,$ is a function from $\,\R\,$ to $\,\R\,$, we have
$\;\R\subset U'\cup V'\,$. 
Thus $\;U'\cap V'\not=\emptyset\;$ because $\,\R\,$ is connected 
and (since $\;U,V\not=\emptyset\,$) both sets $\,U',V'\,$ are nonempty.
Hence we can find $\;a,b,c\in\R\;$ with                
$\;(a,b)\in U\;$ and $\;(a,c)\in V\,$. Since $\;U,V\subset\R^2\;$
are open, there are real numbers $\,u,v\,$ with $\;u<a<v\;$
such that
$\;[u,v]\times\{b\}\,\subset\,U\;$ and $\;[u,v]\times\{c\}\,\subset\,V\,$. 
Now consider the closed set $\;A\,:=\,[u,v]\times\R\setminus(U\cup V)\,$.
We claim that $\;[u,v]\subset \pi(A)\;$ and this yields to
the desired contradiction because if the interval 
$\,[u,v]\,$ is contained in $\,\pi(A)\,$ then  $\;A\cap F\not=\emptyset\;$
contrarily to $\;A\cap F\,\subset\,A\cap(U\cup V)\,=\,\emptyset\,$.
In order to verify $\;[u,v]\subset \pi(A)\;$ let $\;x\in [u,v]\,$.
Then $\;(x,b)\in U\;$ and $\;(x,c)\in V\;$ and hence the connected
set $\,\{x\}\times\R\,$ meets the two disjoint open sets $\,U,V\,$
and therefore $\,\{x\}\times \R\,$ cannot be a subset of $\,U\cup V\,$.
Thus we can find a real number $\,y\,$ such that 
$\;(x,y)\,\in\,(\{x\}\times \R)\setminus(U\cup V)\,\subset\,A\;$
and we arrive at $\;x\in\pi(A)\,$, {\it q.e.d.}
\mp\sp
{\bf Lemma 3.} {\it If $\,Y\,$ is an infinite set 
then there exists a family $\,{\cal Y}\,$ of subsets 
of $\,Y\,$ such that $\,|{\cal Y}|=2^{|Y|}\,$ and 
$\;|A\setminus B|=|Y|\;$ whenever 
$\;A,B\in{\cal Y}\;$ and $\,A\not=B\,$.} 
\mp  
{\it Proof.} Write $\,Y=R\cup S\,$
with $\,|R|=|S|=|Y|\,$ and $\,R\cap S=\emptyset\,$
and let $\,f\,$ be a bijection from $\,R\,$ onto $\,S\,$.
Put
$\;{\cal W}\,:=\,\{\,T\cup(S\setminus f(T))\;|\;T\subset R\,\}\,$.
Obviously, $\,|{\cal W}|=2^{|Y|}\,$ and 
$\;A\not\subset B\;$ whenever 
$\;A,B\in{\cal W}\;$ and $\,A\not=B\,$.
Consequently, the set $\;(A\times Y)\setminus(B\times Y)\,=\,
(A\setminus B)\times Y\;$ is nonempty and hence equipollent with $\,Y\,$
whenever 
$\;A,B\in{\cal W}\;$ and $\,A\not=B\,$.
Therefore, if $\,g\,$ is a bijection from $\;Y\times Y\;$ 
onto $\,Y\,$ then  
$\;{\cal Y}\,:=\,\{\,g(S\times Y)\;|\;S\in{\cal W}\,\}\;$
does the job, {\it q.e.d.}
\mp
{\bf Lemma 4.} {\it If $\,A\,$ is a subspace of $\,\R^2\,$
and $\,\varphi\,$ is a homeomorphism from $\,A\,$ onto 
a subspace of $\,\R^2\,$
then there is a ${\rm G}_\delta$-set $\,G\subset \R^2\,$ 
and a homeomorphism $\,\hat\varphi\,$ from $\,G\,$ onto 
a subspace of $\,\R^2\,$ such that
$\;A\subset G\;$ and 
$\,\hat\varphi(x)=\varphi(x)\,$ for every $\,x\in A\,$.}
\sp
Lemma 4 is an immediate consequence of the 
Lavrentieff expansion theorem [1] 4.3.20. 
\bp
{\bff 4.~Proof of Theorem 1}
\mp
Let $\,{\cal G}\,$ denote the family of all nonempty
G$_\delta$-subsets of $\,\R^2\,$. 
Again, let $\,\pi\,$ be the projection $\;(x,y)\mapsto x\;$
from $\,\R^2\,$ onto $\,\R\,$.
Let $\,{\cal B}\,$ denote the family of all Borel sets $\,B\subset\R^2\,$
such that $\,|\pi(B)|=\bc\,$.
Then $\;|{\cal B}|=|{\cal G}|=\bc\;$ since $\,\R^2\,$ has precisely 
$\,\bc\,$ Borel subsets and 
$\;[0,t]\!\times\!\R\,\in\,{\cal G}\subset{\cal B}\;$
for every $\,t>0\,$.
Let $\,{\cal F}\,$ denote the family of all homeomorphisms 
from sets $\,G\in{\cal G}\,$ onto 
subspaces of $\,\R^2\,$.
(Note that $\,f(G)\in{\cal G}\,$ since $\,X\subset\R^2\,$
is completely metrizable if and only if $\,X\in{\cal G}\,$.)
Since $\,|{\cal G}|=\bc\,$ and since there are only $\,\bc\,$ continuous 
functions from any $\,G\in{\cal G}\,$ into $\,\R^2\,$,
we have $\,|{\cal F}|=\bc\,$.
Define a bijection $\;x\mapsto f_x\;$
from $\,\R\,$ onto $\,{\cal F}\,$
and a bijection $\;x\mapsto B_x\;$
from $\,[0,\infty[\;$ onto $\,{\cal B}\,$.
(The domain $\,[0,\infty[\;$ is chosen for the only reason 
that $\,[0,\infty[\;$ is a set $\,J\subset\R\,$ 
with $\,|J|=|\R\setminus J|=\bc\,$.)
\mp
Let $\,\preceq\,$ be a well-ordering of $\,\R\,$ 
such that $\;|\{\,x\;|\;x\prec y\,\}|<\bc\;$ for every $\,y\in\R\,$.
Via induction we define for every $\,x\in\R\,$
points $\,y_x,z_x\in\R^2\,$ with $\,\pi(y_x)=\pi(z_x)\,$
as follows. 
Assume that for $\,\xi\in\R\,$
points $\,y_x,z_x\in\R^2\,$ with $\,\pi(y_x)=\pi(z_x)\,$ are 
already defined for all $\,x\prec \xi\,$.
Then the definition of $\,y_\xi\,$ and $\,z_\xi\,$
is carried out as follows. 
We distinguish the two cases that either $\,\xi\in[0,\infty[\;$
or $\,\xi\not\in[0,\infty[\,$. (Only for $\,\xi\geq 0\,$
the set $\,B_\xi\,$ is defined!)
\sp
Put $\;Y(\xi)\,:=\,\{\,f_u(y_v)\;|\;u,v\prec\xi\,\}\;$ and 
$\;Z(\xi)\,:=\,\{\,f_u(z_v)\;|\;u,v\prec\xi\,\}\;$ and 
$\;U(\xi)\,:=\,\{\,\pi(y_v)\;|\;v\prec \xi\,\}\,$.
Clearly, the sizes of these three sets are smaller than $\,\bc\,$.
\sp
If $\,\xi\geq 0\,$ then choose a point $\,a\,$ in the nonempty set 
\sp
\cl{$\;B_\xi\setminus\big(Y(\xi)\cup Z(\xi)\cup(U(\xi)\times\R)\big)\;$}
\sp
and define $\,y_\xi=z_\xi=a\,$.
(This set is nonempty since $\;|Y(\xi)\cup Z(\xi)|<\bc\;$
and the size of $\;\pi(U(\xi)\times\R)=U(\xi)\;$ is smaller than $\,\bc\,$,
whereas $\,|\pi(B_\xi)|=\bc\,$.)
\mp
If $\,\xi<0\,$ then 
let $\,\mu\,$ be the minimum of the nonempty set 
$\;\R\setminus U(\xi)\;$
with respect to the well-ordering $\,\preceq\,$
and choose distinct points 
$\,y_\xi,z_\xi\,$ in the vertical line $\,\{\mu\}\times\R\,$
such that $\,\{y_\xi,z_\xi\}\,$ is disjoint 
from $\,Y(\xi)\cup Z(\xi)\,$. 
(This choice is possible since the size of $\;Y(\xi)\cup Z(\xi)\;$
is smaller than $\,\bc\,$ while $\;|\{\mu\}\times\R|=\bc\,$.)
\mp
Now, having defined points $\,y_x,z_x\in\R^2\,$ for every index $\,x\in\R\,$,
put $\;S_y:=\{\,y_x\;|\;x\in\R\,\}\;$ 
and $\;S_z:=\{\,z_x\;|\;x\in\R\,\}\;$ 
and $\;S:=S_y\cup S_z\,$.
By definition, $\,y_u\not=y_v\,$ and $\,z_u\not=z_v\,$  
whenever $\,u\not=v\,$, and $\,y_x\not=z_x\,$ and $\,\pi(y_x)=\pi(z_x)\,$
for every index $\,x<0\,$, and $\,y_x=z_x\,$ for every index $\,x\geq 0\,$.
In particular, $\;|S_y|=|S_z|=|S|=\bc\,$.
Furthermore, $\,\pi(S)=\R\,$ and 
hence $\;1\leq |S\cap (\{x\}\times\R)|\leq 2\;$ 
for every $\,x\in\R\,$.
(It cannot happen by accident that $\,\pi(S)\not=\R\,$
due to the $\preceq$-minimality of the index $\,\mu\,$
and the fact that $\;]{-\infty,0}[\,\not\subset \{\,x\;|\;x\prec y\,\}\;$
for every $\,y\in\R\,$.)
Consequently, for every set $\,T\,$ of negative real indices 
we can define 
a real function $\,F[T]\subset S\,$ by $\;F[T]\,:=\,
\{\,z_x\;|\;x\in T\,\}\cup\{\,y_x\;|\;x\,\in\,\R\setminus T\,\}\,$.
We claim that the following statement is true. 
\mp
(4.1)$\;$ {\it If $\,A,B\subset S\,$
and $\;|A\setminus B|=\bc\;$
then the subspace $\,A\,$ of $\,\R^2\,$ 
cannot be embedded into the subspace $\,B\,$ of $\,\R^2\,$.}
\mp
In order to prove (4.1) let $\,A,B\subset S\,$ and 
let $\,\varphi\,$ be a homeomorphism from $\,A\,$ 
onto a subspace of $\,B\,$. 
By definition, the following two conditions hold.
\sp
(4.2)$\;$ {\it If $\,a,b,x\in\R\,$
and $\,a,b\prec x\,$ 
and $\,y_b\,$ lies in the domain of $\,f_a\,$ 
then $\,y_x\not=f_a(y_b)\,$.} 
\sp
(4.3)$\;$ {\it If $\,a,b,x\in\R\,$
and $\,a,b\prec x\,$ 
and $\,z_b\,$ lies in the domain of $\,f_a\,$ 
then $\,z_x\not=f_a(z_b)\,$.} 
\sp
By Lemma 4 we can expand $\,\varphi\,$
to a homeomorphic embedding $\,f\,$ from $\,G\,$ into $\,\R^2\,$ 
for some $\,G\in{\cal G}\,$ with $\;A\subset G\;$. 
The inverse $\,f^{-1}\,$ is a homeomorphism from $\,f(G)\in{\cal G}\,$
onto $\,G\,$. Thus $\,f,f^{-1}\in{\cal F}\,$ and $\;f(A)\subset B\,$.
\sp
Define $\;\Omega\,:=\,\{\,x\in\R\;|\;y_x\,\in\,A\setminus B\,\}\;$
and choose $\;\gamma,\alpha\in\R\;$ such that
$\,f=f_\gamma\,$ and $\,f^{-1}=f_\alpha\,$ and put 
$\,\Omega_\alpha\,:=\,\{\,x\in\Omega\;|\;x\succ\alpha\,\}\,$.
Let $\,\xi\in\Omega_\alpha\,$ and $\,\beta\in\R\,$ such that
$\;f_\alpha(y_\beta)=y_\xi\,$. Then $\,\beta\succeq\xi\,$ 
in view of (4.2) since $\,\alpha\prec\xi\,$. Moreover, 
$\,\beta\succ \xi\,$ because $\,y_\beta=f_\gamma(y_\xi)\,$ lies in 
$\,B\,$ while $\,y_\xi\not\in B\,$. 
From $\,\xi\prec\beta\,$ and $\,f_\gamma(y_\xi)=y_\beta\,$ and (4.2)
we conclude $\,\beta\preceq\gamma\,$, whence $\,\xi\preceq\gamma\,$.
Thus $\;\Omega_\alpha\subset\{\,\xi\;|\;\xi\preceq\gamma\,\}\;$
and hence $\,|\Omega_\alpha|<\bc\,$.
Trivially, $\;|\Omega\setminus\Omega_\alpha|
\leq|\{\,x\;|\;x\preceq\alpha\,\}|<\bc\,$.
Thus 
$\;|(A\setminus B)\cap S_y|=|\Omega|<\bc\,$.
Similarly, in view of (4.3) we conclude  
$\;|(A\setminus B)\cap S_z|<\bc\,$.
Consequently, from $\,S_x\cup S_y=S\,$ we derive
$\;|A\setminus B|<\bc\;$ and this proves (4.1). 
\mp
Now in order to conclude the proof of Theorem 1,
by applying Lemma 3 we can define 
a family $\,{\cal V}\,$ of subsets 
of $\;]{-\infty,0}[\;$ such that $\,|{\cal V}|=2^{c}\,$ and 
$\;|V\setminus W|=\bc\;$ whenever $\;V,W\in{\cal V}\;$ and $\,V\not=W\,$.
Clearly, if $\,V,W\,$ are sets of negative reals then 
$\;|F[V]\setminus F[W]|=|(V\setminus W)\cup(W\setminus V)|\,$.
Therefore, $\;|F[V]\setminus F[W]|=\bc\;$ 
whenever $\,V,W\in{\cal V}\,$ and $\,V\not= W\,$.
Consequently by virtue of (4.1), for each 
$\,T\in{\cal V}\,$ the space $\,F[T]\,$ cannot be homeomorphic 
to a subspace of $\,F[T']\,$ whenever $\,T\not=T'\in{\cal V}\,$.
Furthermore, for each $\,T\in{\cal V}\,$ the space 
$\,F=F[T]\,$ is not homeomorphic to any 
proper subspace of itself. This is indeed true in view of (4.1)
because if $\,G\,$ is a connected subspace of $\,F\,$ and $\,G\not= F\,$ then 
$\;G\,=\,F\cap(I\times\R)\;$ for some interval $\,I\not=\R\,$
and hence $\;|F\setminus G|=\bc\,$. 
Thus we conclude that the 
$\,2^\bc\,$ spaces $\;F[T]\;(T\in{\cal V})\;$
are incomparable subspaces of $\,\R^2\,$.
\sp
Finally, since $\;F[T]\cap B\not=\emptyset\;$
for every $\,B\in{\cal B}\,$,  
by virtue of Lemma 1 we can be sure that $\,F[T]\,$
intersects every non-null Borel set. A fortiori, 
$\,F[T]\,$ is massive.   
And since every closed subset of $\,\R^2\,$
is a Borel set, for each $\,T\in{\cal V}\,$
the space $\,F[T]\,$ is connected by virtue of Lemma 2.
This concludes the proof of Theorem 1 with the property 
{\it dense} sharpened to {\it massive}.
\mp
{\it Remark.} The spaces $\,F[T]\,$ are not only massive
with respect to the Lebesgue measure on $\,\R^2\,$ 
but also with respect to the Baire property, meaning that 
$\,F[T]\,$ intersects every non-meager subset of $\,\R^2\,$
satisfying the Baire property. (Because, in view of [2] 11.16,
Lemma 1 remains true if the property {\it not null}
is replaced with {\it not meager}.) 
Furthermore, 
all spaces $\,F[T]\,$ are {\it totally imperfect} (see [3]), i.e.~no 
space $\,F[T]\,$ contains an uncountable closed set. 
(Suppose that some uncountable closed subset $\,A\,$ of $\,\R^2\,$ 
lies in $\,F[T]\,$. Then $\,|A|=\bc\,$ and hence $\,|\pi(A)|=\bc\,$ 
since $\,F[T]\,$ is a real function. 
Then $\;B\,=\,\{\,(x,y+1)\;|\;(x,y)\in A\,\}\;$
is a closed set disjoint from $\,F[T]\,$ with $\,\pi(B)=\pi(A)\,$. 
This is impossible since $\,F[T]\,$ meets the Borel set 
$\,B\,$ due to $\,|\pi(B)|=\bc\,$.)
\vfill\eject
\bp
{\bff 5.~Proof of Theorem 2}
\mp
The following lemma is essential for the proofs of Theorems 2 and 3.
\mp
{\bf Lemma 5.} {\it If $\,\tau\in{\cal T}\,$ and the space $\,(\R,\tau)\,$
is connected then every interval is $\tau$-connected.
If $\,F\in{\cal R}\,$ is connected 
then $\;F\cap(I\times\R)\;$ is connected for every interval $\,I\,$.} 
\mp
{\it Proof.} Let $\,\tau\in{\cal T}\,$ and assume that 
the space $\,(\R,\tau)\,$ is connected.
Let $\,a<b\,$. Firstly we verify that $\;I=[a,b]\;$ 
is $\tau$-connected. Assume indirectly that $\,I\,$ is not connected.
Then there are $\tau$-open sets $\;U,V\subset\R\;$ 
such that $\,U\cap I\not=\emptyset\,$ and $\,V\cap I\not=\emptyset\,$
and $\;U\cup V\,\supset\,I\;$ and $\;U\cap V\cap I=\emptyset\,$.
Assume $\,a\in U\,$ and $\,b\in V\,$ (whence $\,a\not\in V\,$
and $\,b\not\in U\,$). 
Then $\;(U\setminus[b,\infty[)\cup{]{-\infty,a}[}\;$
and $\;(V\setminus{]{-\infty, a}]})\cup{]b,\infty[}\;$
are disjoint nonempty $\tau$-open sets which cover the whole set $\,\R\,$. 
This contradicts the connectedness of the space $\,(\R,\tau)\,$.
In a similar way a contradiction can be derived 
if $\,\{a,b\}\,$ is a subset of $\,U\,$ or of $\,V\,$.
So we have proved that 
$\;[a,b]\;$ is $\tau$-connected whenever 
$\,a<b\,$. Since every nondegenerate interval is a union of intervals 
$\;[a_n,b_n]\;(n\in\N)\;$ with $\;a_m\leq a_n<b_n\leq b_m\;$
whenever $\,n<m\,$, we conclude that every interval 
is $\tau$-connected. This proves the first statement of Lemma 5.
The second statement can be verified in a similar way. 
Alternatively, by considering the space $\,(\R,\tau[F])\,$
for $\,F\in{\cal R}\,$, the second statement is a consequence of 
the first, {\it q.e.d.}
\mp
For the proof of Theorem 2 we also need the following noteworthy lemma.
\mp
{\bf Lemma 6.} {\it If $\,F\in{\cal R}\,$ is completely metrizable 
then $\,F\,$ is a nowhere dense subset of $\,\R^2\,$.}
\mp
{\it Proof.} Since $\,\R^2\,$ is a Baire space, it is enough 
to verify that  the closure $\,\overline F\,$
of $\,F\,$ in $\,\R^2\,$ is a meager subset of $\,\R^2\,$.
Since $\,F\,$ is a G$_\delta$-subset of $\,\R^2\,$,
$\,F\,$ is a {\it dense} G$_\delta$-subset of the subspace 
$\,\overline F\,$ of $\,\R^2\,$. 
Hence $\,\overline F\setminus F\,$ is
a {\it meager} F$_\sigma$-subset of the space $\,\overline F\,$.
Therefore, $\,\overline F\setminus F\,$ is a meager subset 
of the plane. (For if $\,A\,$ is a
nowhere dense subset of a closed subspace of $\,\R^2\,$
then $\,A\,$ is clearly nowhere dense in $\,\R^2\,$.) 
Since $\,F\,$ is a G$_\delta$-set 
and $\,|F\cap(\{t\}\times\R)|=1\,$ for every $\,t\in\R\,$,
by the supplementary Kuratowski-Ulam theorem 
[7] 15.4 also $\,F\,$ is a meager subset of $\,\R^2\,$.
Thus, $\,\overline F\,=\,(\overline F\setminus F)\cup F\,$
is a meager subset of $\,\R^2\,$, {\it q.e.d.} 
\mp\sp
Now in order to prove Theorem 2
let $\,F\in{\cal R}\,$ be a connected G$_\delta$-subset
of the Euclidean plane. We have to show that the discontinuity points 
of $\,F\,$ form a meager subset of $\,\R\,$.
Let $\,\pi_2\,$ denote the projection $\;(x,y)\mapsto y\;$
from $\,\R^2\,$ onto $\,\R\,$. 
\sp
Let $\,a\in\R\,$ and assume that 
$\,F\,$ is not continuous at $\,a\in\R\,$.
Then we may choose $\,t_a>0\,$ and 
a sequence $\,x_n\,$ converging to $\,a\,$ 
such that
either $\;F(x_n)> F(a)+t_a\;$ for every $\,n\in\N\,$
or $\;F(x_n)< F(a)-t_a\;$  for every $\,n\in\N\,$.
Put $\;I_n\,:=\,[\min\{a,x_n\}\,,\,\max\{a,x_n\}]\;$
for every $\,n\in\N\,$. Then for every $\,n\in\N\,$ 
the set $\;F\cap(I_n\times\R)\;$
is connected and hence $\;\pi_2(F\cap(I_n\!\times\!\R))\,$
is an interval which contains $\,F(a)\,$ and $\,F(x_n)\,$.
Consequently, either 
$\;\{a\}\times [F(a),F(a)+t_a]\,\subset\,\overline F\;$
or $\;\{a\}\times [F(a)-t_a,F(a)]\,\subset\,\overline F\,$.
\sp
So for every discontinuity point $\,a\,$ of $\,F\,$
there is a nondegenerate closed interval $\,J_a\,$ such that
$\;\{a\}\times J_a\,\subset\,\overline F\,$.
Of course, $\,J_a\,$ is of second category in $\,\R\,$.
By Lemma 6, $\,F\,$ and hence $\,\overline F\,$ is a
nowhere dense subset of $\,\R^2\,$.
Consequently, in view of the Kuratowski-Ulam theorem [7] 15.1,
the set of all discontinuity points of $\,F\,$ must be 
a meager subset of $\,\R\,$.
This concludes the proof of Theorem 2.

\bp
{\bff 6.~Proof of Theorem 3}
\mp
Let $\,\Omega\,$ denote the family of all mappings $\,g\,$ 
from $\,\N\,$ into $\,\{1,2\}\,$ such that $\,g(n)=2\,$
for infinitely many $\,n\in\N\,$.
Naturally, $\,|\Omega|=\bc\,$. 
Transform $\,g\in\Omega\,$ into the double sequence 
\sp
\cl{$\dots,1,1,1,1,1,2,g(1),g(2),g(3),g(4),g(5),\dots$}
\sp
and then replace each comma with the digit $0$ in order to obtain 
a string 
\sp
\cl{$\cdots 1\;0\;1\;0\;1\;0\;1\;0\;1\;0\;2\;
0\;g(1)\;0\;g(2)\;0\;g(3)\;0\;g(4)\;0\;g(5)\;\cdots$}
\sp
which consists of digits from $\,\{0,1,2\}\,$.
Let $\,S(g)\,$ denote this string of digits.
Since in the string $\,S(g)\,$ 
there is clearly a unique first position of the digit $\,2\,$,
we can distinguish the string $\,S(g_1)\,$ from the 
string $\,S(g_2)\,$ whenever 
$\,g_1,g_2\in\Omega\,$ are distinct. 
Moreover, it is clear that every sequence $\,g\in\Omega\,$ 
can be reconstructed from the string $\,S(g)\,$.
Our clue in proving Theorem 3 is 
to define for every $\,g\in\Omega\,$ 
a completely metrizable, connected space $\,F[g]\in{\cal R}\,$ 
such that the string $\,S(g)\,$ (and hence the sequence $\,g\,$)
can be obtained directly from the topology of the space $\,F[g]\,$
by investigating the path-components of $\,F[g]\,$.
There will be three topological types of path-components $\,P\,$
corresponding with the total number $\,n\,$ of {\it noncut points} of $\,P\,$
where $\,n\in\{0,1,2\}\,$. (A point $\,a\,$ in a connected set $\,A\,$
is a noncut point of $\,A\,$ if and only if $\,A\setminus\{a\}\,$ is 
connected.)
\mp
Consider the function 
\sp
\cl{$\;\sigma(u,v)\;:=\;
\{\,(x,\sin{1\over (x-u)(v-x)})\;|\;u<x<v\,\}\;$}
\sp
defined on $\;]u,v[\;$ for $\,u<v\,$.
Clearly, the subspace $\;\{(u,0)\}\cup\sigma(u,v)\cup\{(v,0)\}\;$
of $\,\R^2\,$ is connected and has precisely three path-components, 
namely $\,\{(u,0)\}\,$ and $\,\{(v,0)\}\,$ and $\,\sigma[u,v]\,$.
\mp
For $\,g\in\Omega\,$ define 
$\,F[g]\in{\cal R}\,$ via
\sp
\cl{$F[g]\;\;:=\;\;\{\,(-n,0)\;|\;n\in\N\,\}\,\cup\,([0,1]\times\{0\})\,\cup\,
\bigcup\limits_{k=0}^\infty\sigma(-k-1,-k)
\,\cup\,\bigcup\limits_{k=1}^\infty \Sigma[g,k] $}
\sp
where $\;\Sigma[g,k]\,\subset\;]k,k+1]\times[-1,1]\;$
is defined by 
\sp
\cl{$\;\Sigma[g,k]\,=\,
\sigma(k,k+1)\cup\{(k+1,0)\}\;$} 
\sp
if $\,g(k)=1\,$ and 
\sp
\cl{$\;\Sigma[g,k]\,=\,
\sigma(k,k+{1\over 2})\cup([k+{1\over 2},k+1]\times\{0\})\;$} 
\sp
if $\,g(k)=2\,$. 
\mp 
The set $\,F[g]\,$ is connected because it is clearly possible
to write $\;F[g]\,=\,\bigcup_{n=1}^\infty C_n\;$
where $\,C_n\,$ are connected subsets of $\,\R^2\,$ for all $\,n\in\N\,$
and $\,C_{n+1}\,$ is not separated from $\,C_n\,$ whenever 
$\,n\in\N\,$.                   
With $\,\overline{F[g]}\,$ denoting the closure
of $\,F[g]\,$ in $\,\R^2\,$ we have 
$\;F[g]\,=\,
(S\!\times\!\{0\})\cup(\overline{F[g]}
\setminus(S\!\times\![-1,1]))\;$ for some set 
$\;S\,\subset\,{1\over 2}\Z\;$
and hence $\,F[g]\,$ is a G$_\delta$-subset of $\,\R^2\,$. 
We are done by proving 
that for distinct sequences $\,g,h\in\Omega\,$
the space $\,F[g]\,$ cannot be homeomorphic either to a subspace of $\,F[h]\,$
or to a proper subspace of $\,F[g]\,$.
\mp
For $\,g\in\Omega\,$ let $\,{\cal P}[g]\,$ 
denote the family of all path-components
of $\,F[g]\,$. Let $\,{\cal P}_1[g]\,$ denote the set 
of all {\it singleton} path-components of $\,F[g]\,$.
Thus if $\,P\in{\cal P}_1\,$ then $\;P=\{(x,0)\}\;$ for some $\,x\in\Z\,$.
Let $\,{\cal P}_2[g]\,$ denote the set 
of all {\it infinite and compact} path-components of $\,F[g]\,$
and let $\,{\cal P}_0[g]\,$ denote the set 
of all {\it not compact} path-components of $\,F[g]\,$.
Obviously, if $\,P\in{\cal P}_2\,$ then $\;P\,=\,[a,b]\times\{0\}\;$
for some $\,a<b\,$, and if $\,P\in{\cal P}_0\,$ then 
$\;P\,=\,\sigma(u,v)\;$ for some $\,u<v\,$. 
Consequently, $\,{\cal P}[g]\,$ is the union of the mutually disjoint 
{\it infinite} sets
$\,{\cal P}_0[g]\,$ and $\,{\cal P}_1[g]\,$ and $\,{\cal P}_2[g]\,$.
The notation $\,{\cal P}_n[g]\,$ with $\;n\in\{0,1,2\}\;$ 
corresponds to the fact that 
every space  $\;P\in{\cal P}_n[g]\;$ has precisely $\,n\,$ 
{\it noncut points}. (A singleton has precisely one noncut point
since the empty set is connected vacuously.)
\mp
The family $\,{\cal P}[g]\,$
is naturally ordered
by defining $\,P_1\prec P_2\,$ for $\,P_1,P_2\in{\cal P}[g]\,$
if and only if $\,x_1<x_2\,$ for some $\,(x_1,y_1)\in P_1\,$ 
and some $\,(x_2,y_2)\in P_2\,$.
Obviously the linearly ordered set $\,({\cal P}[g],\prec)\,$ 
is isomorphic with the naturally ordered set $\,\Z\,$.
We claim that the ordering $\,\prec\,$ 
is completely determined by the topology of 
the space $\,F[g]\,$. 
\mp
Obviously there are precisely two 
linear orderings $\,<_1\,$ and $\,<_2\,$ of $\,{\cal P}[g]\,$
such that $\;P,Q\in{\cal P}[g]\;$ are consecutive points with respect 
to $\,<_1\,$ and to $\,<_2\,$ if and only if 
$\,P\not= Q\,$ and $\;P\cup Q\;$ is a connected subset of $\,F[g]\,$. 
One of these linear orderings coincides with $\,\prec\,$
and we can find out which one in a purely topological way. 
Indeed, by definition, 
for every $\,P\in{\cal P}[g]\,$ the set 
$\;\{\,Q\in{\cal P}_2[g]\;|\;Q\prec P\,\}\;$ is finite 
while the set $\;\{\,Q\in{\cal P}_2[g]\;|\;P\prec Q\,\}\;$ is infinite.
This proves the claim.
\mp
Since the ordering $\,\prec\,$ 
is completely determined by the topology of 
the space $\,F[g]\,$, also the sequence 
$\,g\in\Omega\,$ is completely determined 
by the topology of the space $\,F[g]\,$.
To verify this, write the linearly ordered countable set $\,{\cal P}[g]\,$
as a string 
\sp
\cl{$\cdots\;P_{-5}\;P_{-4}\;P_{-3}\;P_{-2}\;P_{-1}\;P_{0}\;
P_{1}\;P_{2}\;P_{3}\;P_{4}\;P_{5}\;\cdots$}
\sp
where $\;P_n\prec P_m\;$ if and only if $\,n<m\,$.
Then replace 
each $\,P_i\,$ with the total number $\,n\,$ of the noncut points
of $\,P_i\,$ in order to obtain a string of digits $\,n\in\{0,1,2\}\,$.
Obviously, this string of digits coincides with the string 
$\,S(g)\,$.
\mp
Due to the previous considerations, 
we can be sure that for distinct $\,g,h\in\Omega\,$ the spaces 
$\,F[g]\,$ and $\,F[h]\,$ cannot be homeomorphic. 
Therefore, the proof is finished by verifying that, whenever $\,g,h\in\Omega\,$,
the space $\,F[h]\,$ cannot be homeomorphic to a {\it proper} subspace
of $\,F[g]\,$.                                       
\mp
Let $\,g,h\in\Omega\,$ and assume 
that $\,\varphi\,$ is a homeomorphism from 
$\,F[h]\,$ onto a proper subspace $\,G\,$
of $\,F[g]\,$.
By Lemma 5, any connected subspace of $\,F[g]\,$ equals 
$\;(I\times\R)\cap F[g]\;$
for some interval $\,I\,$. 
Therefore, $\;G\,=\,(I\times\R)\cap F[g]\;$ for some interval $\,I\not=\R\,$.
Let $\,{\cal P}(G)\,$ denote the family of the path-components 
of the space $\,G\,$. Naturally, 
$\;{\cal P}(G)\,=\,\{\,\varphi(P)\;|\;P\in{\cal P}[h]\,\}\,$.
Let $\,\prec_h\,$ be the linear ordering 
of $\,{\cal P}[h]\,$ and define $\;A\prec_G B\;$ for $\;A,B\in{\cal P}(G)\,$
if and only if $\;\varphi^{-1}(A)\prec_h\varphi^{-1}(B)\,$.
Trivially, the linearly ordered sets $\,({\cal P}[h],\prec_h)\,$ and 
$\,({\cal P}(G),\prec_G)\,$ are isomorphic.
This, however, is not possible. 
Indeed, if $\,<\,$ is any linear ordering 
of $\,{\cal P}(G)\,$ such that 
$\;\bigcup\{\,A\in{\cal P}(G)\;|\;A\leq B\,\}\;$ 
is a connected subspace of $\,G\,$ for every $\,B\in{\cal P}(G)\,$
then the linearly ordered 
set $\,({\cal P}(G),<)\,$ has a maximum or a minimum
since $\,I\not=\R\,$.
On the other hand, the linearly ordered set $\,({\cal P}[h],\prec_h)\,$ 
is isomorphic with $\,\Z\,$ and hence it has neither 
a maximum nor a minimum.
This concludes the proof of Theorem 3.
\mp\sp
{\it Remark.} Concerning a statement
in the abstract, if $\,F_1\,$ and $\,F_3\,$ is any function provided 
by the proof of Theorem 1 and Theorem 3, respectively,  
we can be sure that 
$\,F_1\,$ and $\,F_3\,$ are incomparable.
(Because $\,F_1\,$ is connected and totally pathwise disconnected,
whereas $\,F_3\,$ has infinite path-components and only countably 
many singleton path-components and every connected subset
of a path-component of $\,F_3\,$ is pathwise connected.)
\bp
{\bff 7.~Proof of Theorem 4}
\mp
In the following, an {\it interval partition} is any partition $\,{\cal P}\,$
of the set $\,\R\,$ which consists of intervals. 
(Note that $\,\{a\}=[a,a]\,$ is an interval for every $\,a\in\R\,$.)
For every topology $\,\tau\in{\cal T}\,$ 
let $\,{\cal C}_\tau\,$ denote the family 
of all components of the space $\,(\R,\tau)\,$. 
Furthermore, let $\,{\cal T}^*\,$ denote the family 
of all topologies $\,\tau\in{\cal T}\,$ such that 
$\,(\R,\tau)\,$ is locally connected. 
\sp
First of all note, again, 
that if a topology $\,\tau\,$ on $\,\R\,$
is finer than $\,\eta\,$ then every $\tau$-connected set of reals 
is $\eta$-connected a fortiori. Naturally a set of reals 
is $\eta$-connected if and only if it is an interval.
Consequently, if $\,\tau\in{\cal T}\,$ then the family $\,{\cal C}_\tau\,$
is an interval partition.
\sp
Now assume that $\,\tau\in{\cal T}^*\,$. 
If $\,V\in{\cal C}_\tau\,$ then 
the relative topologies of $\,\tau\,$ and $\,\eta\,$
coincide on the set $\,V\,$ because otherwise, 
by similar arguments as in the proof of Proposition 1,
$\,V\,$ would not be $\tau$-connected. 
Naturally, the elements of $\,{\cal C}_\tau\,$ are 
$\tau$-closed sets. Moreover, since $\,(\R,\tau)\,$ is locally connected,
all members of $\,{\cal C}_\tau\,$ are $\tau$-open as well.
So the space $\,(\R,\tau)\,$ is the topological sum 
of all $\,A\in{\cal C}_\tau\,$, and each $\,A\in{\cal C}_\tau\,$ 
is both a subspace of the space $\,(\R,\tau)\,$ and the Euclidean space 
$\,\R\,$. Consequently, each topology $\,\tau\in{\cal T}^*\,$
is completely determined by the family $\,{\cal C}_\tau\,$ and hence
$\;\tau\mapsto{\cal C}_\tau\;$ 
is an injective mapping from $\,{\cal T}^*\,$ into the 
family of all interval partitions.
\sp
If $\,{\cal P}\,$ is an interval partition then let us call 
a quadruple $\,(\alpha,\beta,\gamma,\delta)\,$
of cardinal numbers the {\it type} of $\,{\cal P}\,$ 
if and only if $\,{\cal P}\,$ contains 
exactly $\,\alpha\,$ 
singletons and exactly $\,\beta\,$ compact intervals
of positive length
and exactly $\,\gamma\,$ half-open intervals
and exactly $\,\delta\,$ open intervals.
Thereby, an interval is {\it half-open} 
when it is of the form $\;]{-v,-u}]\;$ or $\;[u,v[\;$
with $\,u<v\,$ including the possibility $\,v=\infty\,$.
(So we do not care that half-open 
intervals $\,[t,\infty[\;$ and $\;]{-\infty,t}]\;$ are Euclidean closed.)
\sp
It goes without saying that 
for topologies $\,\tau,\sigma\in{\cal T}^*\,$
the spaces $\,(\R,\tau)\,$ and $\,(\R,\sigma)\,$ are 
homeomorphic if and only if the types of the interval partitions 
$\,{\cal C}_\tau\,$ and $\,{\cal C}_\sigma\,$ coincide. 
Consequently, the first statement in Theorem 4 is verified by proving 
that there are only countably many quadruples 
$\,(\alpha,\beta,\gamma,\delta)\,$ which occur as types 
of interval partitions. 
\sp
If $\,(\alpha,\beta,\gamma,\delta)\,$ is the type of an interval partition 
then $\,\beta,\gamma, \delta\,$ are countable cardinals 
or, equivalently, $\;\beta,\gamma, \delta\,\in\,\N\cup\{0,\aleph_0\}\,$.
(Because every nondegenerate interval contains a rational number.) 
Consequently, for every cardinal number $\,\alpha\,$
there are not more than countably many triples 
$\,(\beta,\gamma,\delta)\,$ such that the quadruple 
$\,(\alpha,\beta,\gamma,\delta)\,$ is the type of an interval partition. 
So if there occur only countably many cardinals $\,\alpha\,$
then there occur only countably many quadruples 
$\,(\alpha,\beta,\gamma,\delta)\,$ as possible types. 
For the cardinal number $\,\alpha\,$ we have 
$\,\alpha\leq \bc\,$ since there are precisely $\,\bc\,$ singleton
subsets of $\,\R\,$. Unfortunately 
we cannot rule out that there are uncountably many cardinals below $\,\bc\,$.
Actually, it 
is consistent with ZFC set theory that there are $\,\bc\,$ cardinals below 
$\,\bc\,$. (This can easily be established
by routine forcing.
In G\"odel's universe L 
let $\,\theta\,$ be 
the smallest ordinal number 
of uncountable cofinality satisfying $\,\aleph_\theta=\theta\,$.
While $\,2^{\aleph_0}=\aleph_1\,$ holds in L,
there is a generic extension $\,V[{\rm L}]\,$ of L 
such that
$\,2^{\aleph_0}=\aleph_\theta\,$ holds in $\,V[{\rm L}]\,$, see [2] p.~226.) 
Fortunately, we can prove that if 
$\,(\alpha,\beta,\gamma,\delta)\,$ is the type of an interval partition 
$\,{\cal P}\,$ then $\,\alpha\,$ lies in the countable 
set $\;\N\cup\{0,\aleph_0,\bc\}\,$. 
\sp
Indeed, let $\,{\cal A}\,$ be the family 
of all singletons $\,\{a\}\,$ in $\,{\cal P}\,$, 
whence $\;\alpha=|{\cal A}|\,$. 
Assume that $\,{\cal A}\,$ is uncountable.
Trivially, $\;|{\cal A}|=|\bigcup{\cal A}|\,$. 
Now, $\;\bigcup{\cal A}\;$ is a G$_\delta$-subset 
of $\,\R\,$ because $\;\bigcup{\cal A}\,=\,\R\setminus\bigcup
({\cal P}\setminus{\cal A})\;$ 
and every nondegenerate interval is a F$_\sigma$-subset of $\,\R\,$.
Well, any uncountable G$_\delta$-subset 
of $\,\R\,$ is an uncountable Borel set and hence of size $\,\bc\,$ 
(see [2] 11.20). Consequently, $\;\alpha=|{\cal A}|=\bc\,$.
This concludes the proof of the first statement of Theorem~4.
\mp
Clearly, for $\,\tau\in{\cal T}^*\,$
the locally connected space $\,(\R,\tau)\,$ is separable
if and only if the family 
$\,{\cal C}_\tau\,$ is countable. 
(In other words, $\,\alpha\not=\bc\,$ for the type 
$\,(\alpha,\beta,\gamma,\delta)\,$ of the interval partition  
$\,{\cal C}_\tau\,$.) 
In order to prove the second statement of Theorem 4
let $\,\tau\in{\cal T}^*\,$ and assume that $\,(\R,\tau)\,$ is separable
and let $\,\varphi\,$ be an injective mapping from the countable set
$\,{\cal C}_\tau\,$ into $\,\N\,$. Then define 
$\,F\in{\cal R}\,$ by $\;F(x)=\varphi(V)\;$ for $\,x\in V\,$ 
and $\,V\in{\cal C}_\tau\,$. Clearly, $\,\tau[F]=\tau\,$. 
Furthermore, the set $\,\Gamma(\tau)\,$ is countable
because $\,{\cal C}_\tau\,$ is countable
and the Euclidean interior of any set in the family $\,{\cal C}_\tau\,$
is disjoint from $\,\Gamma(\tau)\,$. 
\sp
To conclude the proof of Theorem 4, we claim that 
$\,\Gamma(\tau)\,$ is a closed subset of $\,\R\,$
if $\,\tau\in{\cal T}^*\,$.
(We do not assume separability and hence we show a bit more than 
stated in Theorem 4.)
\sp
Assume indirectly that 
an Euclidean limit point $\,a\,$ of 
$\,\Gamma(\tau)\,$ does not lie in $\,\Gamma(\tau)\,$.
Then $\;{\cal U}_\eta(a)={\cal U}_\tau(a)\;$ 
and $\,{\cal U}_\tau(a)\,$ has a basis $\,{\cal B}(a)\,$
consisting of $\tau$-connected and Euclidean open neighborhoods of $\,a\,$.
Let $\,B\in{\cal B}(a)\,$. Then $\;B=]u,v[\;$ for $\,u<a<v\,$
and $\;\Gamma(\tau)\cap B\not=\emptyset\,$. 
Let $\;x\,\in\,\Gamma(\tau)\cap B\,$.
Then we can find a $\tau$-connected, $\tau$-open neighborhood $\,V\,$
of $\,x\,$ such that $\;x\in V\subset B\;$ 
and $\,V\,$ is not $\eta$-open,
whence $\,V=\{x\}\,$
or $\;V=[a,b[\;$ or $\;V=\,]a,b]\;$ or $\;V=[a,b]\;$ 
for $\,a<b\,$. This enables us to write $\,B\,$
as a disjoint union of nonempty $\tau$-open subsets 
and hence we arrive at the contradiction 
that $\,B\,$ is not $\tau$-connected.
\sp
This finishes the proof of Theorem 4.
Furthermore, since any interval is clearly a 
locally compact and completely metrizable subspace of the real line,
we additionally have proved the following noteworthy 
proposition.
\mp
{\bf Proposition 2.} {\it 
If $\,\tau\in{\cal T}\,$ and $\,(\R,\tau)\,$ is locally connected
then (not assuming separability) 
$\,(\R,\tau)\,$ is locally compact and completely metrizable
and the set $\,\Gamma(\tau)\,$ is Euclidean closed. 
(And $\,\Gamma(\tau)\,$ is countable if and only if 
the locally compact space $\,(\R,\tau)\,$ is second countable.)}
\mp
{\it Remark.} There exists a topology $\,\tau\in{\cal T}\,$
such that the space $\,(\R,\tau)\,$ is locally compact and separable
but neither locally connected nor metrizable, see [8] No.~65.
If $\,F\in{\cal R}\,$ and $\,a\in F\,$
and the space $\,F\,$ has a local basis of connected sets at 
$\,a\,$ then, despite Proposition 2, one cannot conclude 
that the point $\,a\,$ has a compact neighborhood.
For a counterexample consider the completely metrizable, 
connected space $\,F\in{\cal R}\,$
defined in the following way.  With $\,\sigma(u,v)\,$ as in Section 6
put $\;g_n=\sigma(2^{-n},2^{-n+1})\;$ for every $\,n\in\N\,$ 
and define $\,F\in{\cal R}\,$ by $\;F(x)=2^{-n}g_n(x)\;$
if $\,n\in\N\,$ and $\,2^{-n}<x<2^{-n+1}\,$ and $\,F(x)=0\,$ otherwise.
\bp
{\bff 8.~Partitions of the real line}
\mp
In a very natural way the considerations in the proof of Theorem 4 lead to 
a complete classification of all
locally connected refinements of the real line.
First of all we point out that the injective mapping 
$\;\tau\mapsto {\cal C}_\tau\;$ 
from $\,{\cal T}^*\,$ into the family of all interval partitions
is bijective. Because if $\,{\cal P}\,$ is any interval partition 
then there is a unique topology $\,\tau\,$ in the 
family $\,{\cal T}^*\,$ with $\,{\cal C}_\tau={\cal P}\,$.  
This topology $\,\tau\,$ is generated by the basis $\,{\cal B}\,$
where $\,B\in{\cal B}\,$ if and only if 
$\;B=J\cap A\;$ for some open interval $\,J\,$
and some $\,A\in{\cal P}\,$.
(Alternatively, $\,\tau\,$ is generated by the metric $\,d\,$
defined by $\;d(x,y)\,=\,\arctan|x-y|\;$ when $\,x,y\,$ lie in a common 
set $\,A\in{\cal P}\,$ and $\,d(x,y)=2\,$ otherwise.)
Since the mapping $\;\tau\mapsto {\cal C}_\tau\;$ is bijective
and since each topology $\,\tau\in{\cal T}^*\,$ 
is completely determined up to homeomorphism 
by the type of the interval partition $\,{\cal C}_\tau\,$,
in order to establish a classification of all 
topologies $\,\tau\in{\cal T}^*\,$
we have to find out which quadruples 
of cardinal numbers actually occur as types of interval partitions.
It is reasonable to distinguish beteen countable and 
uncountable interval partitions because for $\,\tau\in{\cal T}^*\,$ 
the space $\,(\R,\tau)\,$ is separable if and only if 
the interval partition $\,{\cal C}_\tau\,$ is countable.
\sp
Let  $\;\Omega_1\,:=\,\N\cup\{0,\aleph_0\}\;$ denote the set of 
all countable cardinal numbers. 
Note that $\,\aleph_0+\aleph_0=\aleph_0\,$ 
and $\,n+\aleph_0=n\cdot\aleph_0=\aleph_0\,$ for every $\,n\in\N\,$.
Define sets of quadruples of cardinal numbers by  
$\;{\cal Q}_2\,:=\,\{\bc\}\times\Omega_1^3\;$ and 
\sp
\cl{$\;{\cal Q}_1\;:=\;\big\{\,(\alpha,\beta,\gamma,\delta)\in\Omega_1^4\;\;
\big|\;\;\alpha+\beta+1=\delta\;\lor\;
(\alpha+\beta+1>\delta\;\land\;\gamma=\aleph_0)\,\big\}\,$.}
\mp\sp
{\bf Theorem 5.} {\it If $\,{\cal P}\,$ is an interval partition 
of type $\,\vartheta\,$ then $\,\vartheta\in{\cal Q}_1\,$
when $\,{\cal P}\,$ is countable and $\,\vartheta\in{\cal Q}_2\,$
when $\,{\cal P}\,$ is uncountable.
For every $\;\vartheta\,\in\,{\cal Q}_1\cup{\cal Q}_2\;$
there exists an interval partition of type $\,\vartheta\,$.}
\mp
In order to prove Theorem 5 let 
$\,{\cal P}\,$ be an interval partition of type 
$\,(\alpha,\beta,\gamma,\delta)\,$.
As already verified, 
$\;\beta,\gamma,\delta\,\in\,\Omega_1\;$ 
and $\,\alpha=\bc\,$ 
if and only if $\,{\cal P}\,$ is uncountable.
In particular, if $\,{\cal P}\,$ is uncountable
then $\,(\alpha,\beta,\gamma,\delta)\in{\cal Q}_2\,$.
Conversely, if $\,(\alpha,\beta,\gamma,\delta)\in{\cal Q}_2\,$
(and hence $\,\alpha=\bc\,$)
then we can easily find an (uncountable) interval partition of type 
$\,(\alpha,\beta,\gamma,\delta)\,$. 
Indeed, let $\,\beta,\gamma,\delta\,$ be countable cardinals   
and put $\;{\cal P}\,=\,{\cal P}_0\cup
{\cal P}_1\cup{\cal P}_2\cup{\cal P}_3\;$
where 
\sp
\qquad $\;{\cal P}_1\;=\;\{\,[6n,6n\!+\!1]\;|\;n\in\Z\;\land\;
0\leq n<\beta\,\}\;$ 
\quad and 
\sp
\qquad $\;{\cal P}_2\;=\;\{\,[6n\!+\!2,6n\!+\!3[\,\;|\;
n\in\Z\;\land\;0\leq n<\gamma\,\}\;$ 
\quad and
\sp
\qquad $\;{\cal P}_3\;=\;
\{\,]6n\!+\!4,6n\!+\!5[\,\;|\;n\in\Z\;\land\;0\leq n<\delta\,\}\;$ 
\quad and 
\sp
\qquad $\;{\cal P}_0\;=\;\{\,\{x\}\;\,|\;\,x\,\in\,\R\setminus
\bigcup({\cal P}_1\cup{\cal P}_2\cup{\cal P}_3)\,\}\;$.
\mp
Thus the proof of Theorem 5 is finished by verifying the following two 
statements.
\mp
(8.1)\quad{\it Every quadruple in $\,{\cal Q}_1\,$
is the type of some countable interval partition.}
\sp
(8.2)\quad{\it The type of each countable 
interval partition lies in $\,{\cal Q}_1\,$.}
\mp
In the following, $\,\alpha,\beta,\gamma,\delta\,$ 
are always countable cardinals
and $\,(\alpha,\beta,\gamma,\delta)\,$ 
is {\it admissible} if and only if $\,(\alpha,\beta,\gamma,\delta)\,$ is 
the type of some countable interval partition. 
First of all we point out that if $\;\alpha+\beta+1=\delta\;$
then $\,(\alpha,\beta,0,\delta)\,$ is admissible. 
This is true because 
\sp
\cl{$\;W\;:=\;\bigcup\limits_{0\leq n<\alpha}\!\!\{3n\}\;\,\cup
\bigcup\limits_{0\leq n<\beta}\!\![3n\!+\!1,3n\!+\!2]\;$}
\sp
is Euclidean closed and the open set $\;\R\setminus W\;$
is a union of precisly $\,\delta\,$ mutually disjoint nonempty 
open intervals. 
\mp
If $\,J\,$ is an infinite interval 
and $\,\kappa\,$ is a cardinal number with $\,\kappa\leq\aleph_0\,$
then we clearly can find a subinterval $\,\tilde J\,$
of $\,J\,$ homeomorphic with $\,J\,$
such that $\;J\setminus\tilde J\;$ can be written  as a disjoint union 
of precisely $\,\kappa\,$ half-open intervals.
In case that $\,\kappa=\aleph_0\,$ and $\,J\,$ is open, 
the interval $\,J\,$ itself is such a disjoint union. 
Therefore and since it is clear that $\,\beta+\gamma+\delta>0\,$,
the following statement is true.
\sp
(8.3)\quad {\it If $\,(\alpha,\beta,\gamma,\delta)\,$ is admissible
then $\,(\alpha,\beta,\gamma',\delta')\,$ is admissible 
for all pairs $\,(\gamma',\delta')\,$ of countable cardinals satisfying either
$\,\gamma'\geq \gamma\,$ and $\,\delta'=\delta\,$,
or $\,\gamma'=\aleph_0\,$ and $\,\delta'<\delta\,$.}
\mp
Now in order to prove (8.1) 
assume that $\,\vartheta=(\alpha,\beta,\gamma,\delta)\,$
lies in $\,{\cal Q}_1\,$.
We already know that $\,\vartheta\,$ is admissible if 
$\;\alpha+\beta+1=\delta\,$ and $\,\gamma=0\,$.
Thus by (8.3) $\,\vartheta\,$ is admissible if 
$\;\alpha+\beta+1=\delta\,$ and $\,\gamma\,$ is arbitrary.
In particular, $\,(\alpha,\beta,\aleph_0,\alpha+\beta+1)\,$
is admissible. Thus by (8.3) $\,\vartheta\,$ is admissible
whenever $\;\alpha+\beta+1>\delta\,$ and $\,\gamma=\aleph_0\,$.
This concludes the proof of (8.1).
\mp
Before we prove (8.2) we verify the following two statements.
\mp
(8.4)\quad {\it If $\,(\alpha,\beta,\gamma,\delta)\,$ is admissible
then $\,(\alpha+2\beta+\gamma,0,0,\beta+\gamma+\delta)\,$ is admissible.} 
\sp
(8.5)\quad {\it If $\,(\alpha,\beta,\gamma,\delta)\,$ is admissible
and $\,\alpha,\beta,\gamma\,$ are finite 
then $\,\delta=\alpha+\beta+1\,$.} 
\mp
To verify (8.4) let $\,{\cal P}\,$ be an interval partition 
of type $\,(\alpha,\beta,\gamma,\delta)\,$ and define $\,B\subset\R\,$
such that $\,x\in B\,$ if and only if $\,x\,$ is a noncut point 
of the Euclidean space $\,I\,$ 
for some $\,I\in{\cal P}\,$.
Clearly, $\,|B|=\alpha+2\beta+\gamma\,$. 
Then $\;(\{\,I\setminus B\;|\;I\in{\cal P}\,\}\setminus\{\emptyset\})
\cup\{\,\{b\}\;|\;b\in B\,\}\;$
is an interval partition which consists of precisely $\,|B|\,$ 
singletons and precisely $\,\beta+\gamma+\delta\,$ open intervals.
\sp
To verify (8.5) let $\,(\alpha,\beta,\gamma,\delta)\,$ be admissible.
Clearly, if $\,E\subset\R\,$ is finite then 
$\;\R\setminus E\;$ is a disjoint union of precisely $\,|E|+1\,$ open 
intervals. Therefore and by applying (8.4) we derive 
$\;\alpha+2\beta+\gamma+1=\beta+\gamma+\delta\;$
and hence $\,\delta=\alpha+\beta+1\,$ if
$\,\alpha,\beta,\gamma\,$ are finite. 
\mp
Now in order to prove (8.2) assume that 
$\,\vartheta=(\alpha,\beta,\gamma,\delta)\,$
is admissable. 
By definition, $\,\vartheta\in{\cal Q}_1\,$ 
if $\,\alpha+\beta+1=\delta\,$. In the following  
we will check the two remaining cases 
$\,\alpha+\beta+1<\delta\,$ and $\,\alpha+\beta+1>\delta\,$. 
\sp
Assume firstly that $\,\alpha+\beta+1>\delta\,$. 
Then $\,\alpha+\beta+\gamma=\aleph_0\,$ 
because $\,\alpha+\beta+\gamma<\aleph_0\,$ 
implies $\,\alpha+\beta+1=\delta\,$ by (8.5).
If $\,\gamma<\aleph_0\,$ then $\,\alpha+\beta=\aleph_0\,$
(since $\,\alpha+\beta+\gamma=\aleph_0\,$) 
and $\,\delta<\aleph_0\,$ (since $\,\alpha+\beta+1>\delta\,$)
and hence $\,\gamma+\delta<\aleph_0\,$.
But $\,\alpha+\beta=\aleph_0>\gamma+\delta\,$
is not possible because 
if $\,K_1,K_2,K_3...\,$ is a sequence of mutually
disjoint nonempty compact intervals (of positive length or of length 0)
then $\;\R\setminus \bigcup_{n=1}^\infty K_n\;$ is nonempty and 
cannot be written as a union of finitely many intervals.
(This fact is not so evident as it might seem at first sight. 
To prove this fact note that, by applying  
Sierpi\'nski's theorem [1] 6.1.27,
if $\,u\in K_1\,$ and $\,v\in K_2\,$ and $\,u<v\,$
then $\;[u,v]\not=[u,v]\cap\bigcup_{n=1}^\infty K_n\,$.)
So we can be sure that if $\,\alpha+\beta+1>\delta\,$
then $\,\gamma=\aleph_0\,$ and hence $\,\vartheta\in{\cal Q}_1\,$.
Furthermore, (8.2) is settled by verifying 
that the case $\,\alpha+\beta+1<\delta\,$ can be ruled out. 
\sp
In order to verify this, assume $\,\alpha+\beta+1<\delta\,$.
Then $\;\alpha+\beta=m\;$
for some $\,m\in\N\,$. Let $\,{\cal P}\,$ be an interval partition 
of type $\,\vartheta\,$ and let $\,K\,$ be the union 
of the $\,m\,$ compact elements of $\,{\cal P}\,$.
Then $\,K\,$ is a compact subspace 
of $\,\R\,$ with precisely $\,m\,$ components. 
Let $\,{\cal U}\,$ be the family of the components 
of the open subspace $\,\R\setminus K\,$ of $\,\R\,$. 
Naturally, each set $\,U\in{\cal U}\,$ is a nonempty open interval 
and $\,|{\cal U}|=m+1\,$.
Furthermore, each interval $\,U\in{\cal U}\,$ 
is a disjoint union of non-compact elements of $\,{\cal P}\,$.
Therefore, from $\;|{\cal U}|=\alpha+\beta+1<\delta\;$ 
we conclude that some open interval $\,V\in{\cal U}\,$
must contain two disjoint {\it open} intervals $\,I_1,I_2\in{\cal P}\,$
with $\;\sup I_1=a\leq b=\inf I_2\,$. 
Then $\,[a,b]\subset V\,$ and
$\,[a,b]\,$ is disjoint from $\,K\,$.
In particular, $\,a<b\,$. If $\,{\cal J}\,$ is the family of
all {\it open} $\,J\in{\cal P}\,$ with $\,J\subset [a,b]\,$ 
then we further conclude that each component $\,[u,v]\,$ of the 
non-empty compact set $\;[a,b]\setminus\bigcup{\cal J}\;$ 
is a disjoint union of half-open intervals.
But it is straightforward to verify 
that a compact interval cannot be decomposed
into half-open intervals. This finishes 
the proof of (8.2) and of Theorem 5.
\bp
The classification of the locally connected refinements of the real 
line leads to the following interesting proposition 
about incomparability and local connectedness.
\mp
{\bf Proposition 3.} {\it Let $\,\tau_1,\tau_2\in{\cal T}\,$ be 
such that $\,(\R,\tau_i)\,$ is locally connected for $\,i\in\{1,2\}\,$.
If $\,(\R,\tau_1)\,$ is separable and $\,(\R,\tau_2)\,$ is not discrete
then $\,(\R,\tau_1)\,$ is homemorphic to a subspace 
of $\,(\R,\tau_2)\,$. If $\,(\R,\tau_2)\,$ is not separable and not discrete
then $\,(\R,\tau_1)\,$ is homemorphic to a subspace 
of $\,(\R,\tau_2)\,$.}
\mp
{\it Proof.} Let $\,\tau\in{\cal T}\,$ such that 
$\,(\R,\tau)\,$ is locally connected and separable  
and $\,(\alpha,\beta,\gamma,\delta)\,$ is the type of 
the family $\,{\cal C}_\tau\,$ of all components of the space 
$\,(\R,\tau)\,$.
With $\;n\,\in\,\N\cup\{0\}\;$
define a subspace $\,X\,$ of the real line by 
\sp
\cl{$X\;:=\;\bigcup\limits_{n<\alpha}\{7n\}\;\cup\;
\bigcup\limits_{n<\beta}[7n+1,7n+2]\;\cup\;
\bigcup\limits_{n<\gamma}[7n+3,7n+4[\,\;\cup\;
\bigcup\limits_{n<\delta}]7n+5,7n+6[\,.$}
\sp
It goes without saying that the spaces $\,(\R,\tau)\,$ and 
$\,X\,$ are homeomorphic. If $\,\tau_2\in{\cal T}\,$ such that 
$\,(\R,\tau_2)\,$ is locally connected and not discrete
then some component $\,C\,$ of the space is not a singleton 
and hence the real line $\,\R\,$ is homeomorphic to a subspace
of $\,C\,$. Consequently, by transitivity, 
$\,(\R,\tau)\,$ is homeomorphic to a subspace of $\,(\R,\tau_2)\,$.
\sp
Now let $\,\tau_2\in{\cal T}\,$ such that 
$\,(\R,\tau_2)\,$ is not separable and not discrete
then $\,(\R,\tau_2)\,$ has an open-closed discrete subset $\,D\,$
of size $\,\bc\,$ and the $\tau_2$-open-closed set $\,\R\setminus D\,$ 
contains a copy of the real line. 
If $\,\tau_1\in{\cal T}\,$ such that 
$\,(\R,\tau_1)\,$ is separable then we already know that 
$\,(\R,\tau_1)\,$ is homeomorphic to a subspace of the real line 
and hence $\,(\R,\tau_1)\,$ is homeomorphic to a subspace of 
the subspace $\,\R\setminus D\,$ of $\,(\R,\tau_2)\,$.
If $\,\tau_1\in{\cal T}\,$ such that 
$\,(\R,\tau_1)\,$ is not separable then 
$\,(\R,\tau_1)\,$ is the topological sum of a discrete space of 
size $\,\bc\,$ and a locally connected subspace of $\,\R\,$.
Hence $\,(\R,\tau_1)\,$ is homeomorphic to a subspace of the 
space $\,(\R,\tau_2)\,$, {\it q.e.d.}
\bp
{\bff 9.~An application of Theorem 1}
\mp
A trivial consequence of Theorem 1 is the following corollary 
whose direct proof would not be significantly easier than our
proof of Theorem 1. 
\sp
{\bf Corollary 3.} {\it There exists a real function $\,F\,$
such that $\,F\,$ is a connected, dense subspace of $\,\R^2\,$
and $\,F\,$ is not homeomorphic to a proper subspace of $\,F\,$.}
\mp
In this final section we apply Corollary 3 
in order to prove a noteworthy theorem 
about separable and connected 
refinements of the real line which 
are not of the form $\,(\R,\tau[F])\,$ for a real function 
$\,F\,$. In doing so we realize that  
in Corollary 1 the cardinality $\,2^\bc\,$ can be 
increased to the greatest possible cardinality 
if the property {\it metrizable} is omitted.
\mp
{\bf Theorem 6.} {\it There exist $\,2^{2^\bc}\,$
topologies $\,\tau\in{\cal T}\,$ with $\,\Gamma(\tau)=\R\,$
such that the corresponding Hausdorff spaces $\,(\R,\tau)\,$ 
are incomparable, separable and connected.}
\mp
{\it Proof.} In the following, call a set $\,S\,$ of reals $\bc$-dense
if and only if $\,|S\cap[u,v]|=\bc\,$ whenever $\,u<v\,$.
For the sake of simply writing we use terminology from linear algebra
and regard $\,\R\,$ as a vector space over the field $\,\Q\,$.
Fix a set $\,Y\,$ of irrational numbers such that 
$\,\{1\}\cup Y\,$ is a basis of the vector space $\,\R\,$ over $\,\Q\,$.
Naturally, $\,|Y|=\bc\,$.
For $\,A\subset Y\,$ let $\,L(A)\,$ denote 
the linear subspace of $\,\R\,$ 
generated by the vectors in the set 
$\,\{1\}\cup A\,$. 
Obviously, $\;\Q\subset L(A)\;$ for every $\,A\subset Y\,$
and if $\,A_1,A_2\subset Y\,$ are disjoint then 
$\,L(A_1)\cap L(A_2)=\Q\,$. 
\sp
Let $\,{\bf U}\,$ denote the family of all ultrafilters $\,{\cal U}\,$ 
on $\,Y\,$ such that $\,|A|=\bc\,$ for every $\,A\in{\cal U}\,$.
Then $\,|{\bf U}|=2^{2^{\bc}}\,$ by [2] Theorem 7.6.
It is plain that $\,L(A)\,$ is $\bc$-dense
whenever  $\,A\subset Y\,$ and $\,|A|=\bc\,$.
In particular, if $\,{\cal U}\in{\bf U}\,$ 
then $\,L(A)\,$ is $\bc$-dense for every $\,A\in{\cal U}\,$.
We claim that the following is true for every $\,{\cal U}\in{\bf U}\,$.
\sp
(9.1)$\;$ {\it For every $\,x\in\R\,$ there exists a set $\,A\in{\cal U}\,$
such that $\,x\in L(A)\,$ and $\,[u,v]\setminus L(A)\,$ 
is an infinite set whenever $\,u<v\,$.}
\sp
In order to prove (9.1)
fix a countably infinite set $\,Y_0\subset Y\,$. Then 
$\,Y_0\not\in{\cal U}\,$ and hence $\;Y\setminus Y_0\,\in\,{\cal U}\,$.
Let $\,x\in \R\,$. Then  $\,x\in L(T)\,$ 
for some finite $\,T\subset Y\,$.
Put $\;A\,=\,(Y\setminus Y_0)\cup T\,$. Then 
$\,A\in{\cal U}\,$ and $\,x\in L(A)\,$.
Furthermore, $\;Y\setminus A\not=\emptyset\;$ since $\,|T|<\aleph_0=|Y_0|\,$.  
Thus the set $\;L(Y\setminus A)\setminus\Q\;$ is a dense subset of 
$\,\R\,$. And, clearly,                          
$\;L(Y\setminus A)\setminus\Q\,\subset\,\R\setminus L(A)\,$. 
\mp
Put $\;J:=[-1,1]\,$ and 
let $\,\eta(J)\,=\,\{\,U\cap J\;|\;U\in\eta\,\}\;$ 
denote the relative topology of $\,\eta\,$ on the interval $\,J\,$.
For every ultrafilter $\,{\cal U}\in{\bf U}\,$
let $\,\eta[{\cal U}]\,$ denote the topology on the set $\,J\,$
generated by the subbasis 
$\;\eta(J)\,\cup\,\{\,L(A)\cap J\;|\;A\in{\cal U}\,\}\,$.
\sp
Put $\,\tau=\eta[{\cal U}]\,$. 
Since $\,L(A)\,$ is $\bc$-dense and $\,\Q\subset L(A)\,$
for every $\,A\in{\cal U}\,$,
it is clear that every nonempty open subset of the space 
$\,(J,\tau)\,$ is a set of size $\,\bc\,$
which intersects the countable set $\,\Q\cap J\,$. Consequently, the 
space $\,(J,\tau)\,$ is separable.
From (9.1) we conclude that 
every $\,x\in J\,$ has a $\tau$-neighborhood which is not 
an $\eta(J)$-neighborhood. Thus $\,\tau\,$ is strictly finer than 
$\,\eta(J)\,$ at every point in $\,J\,$.
A moment's reflection suffices to see
that every subinterval of $\,J\,$ is connected in 
the space $\,(J,\tau)\,$. In view of (9.1) it is straightforward that
a nondegenerate subinterval of $\,J\,$
is never a regular subspace of $\,(J,\tau)\,$.
\sp
Assume that $\,{\cal U}_1,{\cal U}_2\in{\bf U}\,$ are
distinct. Then we can choose $\,A\subset\R\,$ such that
$\,A\,$ lies in $\,{\cal U}_1\,$ and 
$\,Y\setminus A\,$ lies in $\,{\cal U}_2\,$. Hence
$\;W_1\,:=\,J\cap L(A)\;$ is $\eta[{\cal U}_1]$-open  and
$\;W_2\,=\,J\cap L(Y\setminus A)\;$ is 
$\eta[{\cal U}_2]$-open. Therefore, the topologies  
$\,\eta[{\cal U}_1]\,$ and $\,\eta[{\cal U}_2]\,$ are distinct as well
because 
$\;\eta[{\cal U}_1]=\eta[{\cal U}_2]=\tau\;$
implies that $\,W_1\,$ and $\,W_2\,$ are $\tau$-open 
sets and hence $\;W_1\cap W_2\,=\,J\cap\Q\;$ is a nonempty 
$\tau$-open set of size smaller than $\,\bc\,$
whose existence has already been ruled out.
Consequently, the family 
\sp
\cl{$\;{\cal G}\;:=\;
\{\,\eta[{\cal U}]\;\;|\;\;{\cal U}\in{\bf U}\,\}\;$}
\sp
consists of $\,2^{2^{\bc}}\,$ topologies $\,\tau\,$ on $\,J=[-1,1]\,$
everywhere strictly finer than $\,\eta(J)\,$ 
such that the spaces $\,(J,\tau)\,$ are separable 
and connected and not regular but, a fortiori, Hausdorff.
\sp
Define an equivalence relation $\,\sim\,$ on $\,{\cal G}\,$
such that $\,\tau_1\sim\tau_2\,$ if and only if 
the spaces $\,([-1,1],\tau_1)\,$ and $\,([-1,1],\tau_2)\,$ are 
homeomorphic. Since there are only $\,2^{\bc}\,$ permutations 
of the set $\,[-1,1]\,$, the size of an equivalence class
cannot exceed $\,2^\bc\,$. Therefore,
since $\;|{\cal G}|=2^{2^\bc}>2^\bc\,$, 
there are $\,2^{2^\bc}\,$ equivalence classes. 
Hence by choosing one topology in each equivalence class
we can select a family  
$\;{\cal E}\subset{\cal G}\;$ such that 
$\;|{\cal E}|=|{\cal G}|=2^{2^\bc}\;$ 
and all spaces $\;([-1,1],\tau)\;(\tau\in{\cal E})\;$
are {\it mutually non-homeomorphic.}
\mp
Now we are ready to verify Theorem 6.
Put $\;I_{-1}:=\,]{-\infty,-1}[\;$ 
and $\;I_1:=\,]1,\infty[\,$. 
Let $\,F\,$ be a real function as depicted in Corollary~3
and consider the separable, connected, metrizable 
topology $\,\tau[F]\in{\cal T}\,$.
Since both intervals can easily be mapped 
onto $\,\R\,$ by strictly increasing functions, 
with the help of the one topology 
$\,\tau[F]\,$ on $\,\R\,$ 
for both indices $\,k=\pm 1\,$ 
we can define a topology $\,\mu_k\,$ on $\,I_k\,$ 
strictly finer than the restriction of $\,\eta\,$ to 
$\,I_k\,$ such that $\,(I_k,\mu_k)\,$ 
is homeomorphic to the space $\,F\,$.
For every topology $\,\tau\in{\cal E}\,$ define a 
separable Hausdorff topology 
$\,\hat\tau\,$ on the set $\,\R\,$ in the following way. 
Firstly, for $\,k\in\{-1,1\}\,$ 
the subspace $\,I_k\,$ of $\,(\R,\hat\tau)\,$ coincides with the space
$\,(I_k,\mu_k)\,$ and the subspace 
$\;[{-1,1}]\;$ of $\,(\R,\hat\tau)\,$ coincides with the 
space $\,([-1,1],\tau)\,$.
Secondly, for $\,k\in\{-1,1\}\,$ 
the family $\;\{\,]k-t,k+t[\cap V\;|\;t>0\,,\;k\in V\in\tau\}\;$
is a local basis at the point $\,k\,$.
Obviously, $\,\hat\tau\,$
is strictly finer than $\,\eta\,$ at every $\,x\in\R\,$.
In view of Lemma 5 it is clear 
that the separable Hausdorff space $\,(\R,\hat\tau)\,$ is connected 
and, moreover, the connected sets are precisely the intervals. 
\sp
Let $\,\tau_1,\tau_2\,$ be topologies in $\,{\cal E}\,$
and assume that $\,f\,$ is a homeomorphism from the connected space 
$\;(\R,\hat\tau_1)\;$ onto a subspace $\,S\,$ of $\;(\R,\hat\tau_2)\,$.
Then $\,f\,$ maps every interval onto some interval and hence
$\,f\,$ is a {\it monotonic} bijection from $\,\R\,$ onto 
the interval $\,S\,$.
Therefore, since $\,x\in\R\,$ has a metrizable neighborhood
in the space $\,(\R,\hat\tau_i)\,$ if and only if 
$\,x\not\in[-1,1]\,$ and since the subspace 
$\,S\setminus f([-1,1])\,$ of $\,(\R,\hat\tau_2)\,$
has precisely two components
and both components are metrizable spaces,
we conclude that $\,f([-1,1])=[-1,1]\,$. This implies 
firstly $\,\tau_1=\tau_2\,$ and, secondly,
$\;S\cap I_k=I_k\;$ for $\,k=\pm1\,$
in view of Corollary 3
and since the two spaces $\,(I_{\pm 1},\mu_{\pm 1})\,$ 
are homeomorphic with $\,F\,$. Hence 
we obtain $\,\hat\tau_1=\hat\tau_2\,$ and $\,S=\R\,$
and arrive at the conclusion that 
the $\,2^{2^\bc}\,$ spaces $\;(\R,\hat\tau)\;(\tau\in{\cal E})\;$
are incomparable, {\it q.e.d.}
\vfill\eject
\bp\bp\bp
{\bff References}
\mp
[1] Engelking, R.: General Topology, revised and completed edition.
Heldermann, 1989. 
\sp
[2] Jech, T.: Set Theory. 3rd ed. Springer 2002.
\sp
[3] Kharazishvili, A.B.: Strange Functions in Real Analysis.
Chapman and Hall 2006.
\sp
[4] Kuba, G.: {\it Counting ultrametric spaces.}
Colloq.~Math. {\bf 152} (2018), 217-234.
\sp
[5] Kulpa, W.: {\it On the existence of maps having graphs connected 
and dense.}

\rl{Fund.~Math.~{\bf 76} (1972), 207-211.}
\sp
[6] Kuratowski, K.: Topology, vol.~1. Academic Press 1966. 
\sp
[7] Oxtoby, J.C.: Measure and Category. Springer 1971. 
\sp

[8] Steen, L.A., and Seebach Jr., J.A.,
Counterexamples in Topology. Dover 1995.
\bp\bp\bp
{\sl Author's address:} Institute of Mathematics,

University of Natural Resources and Life Sciences, Vienna, Austria. 
\sp
{\sl E-mail:} {\tt gerald.kuba(at)boku.ac.at}
\bp\bp\bp
\hrule
\bp\bp\bp
{\bf This paper will be published in ANALYSIS MATHEMATICA.}
\end

\end